\documentclass{article}

\usepackage{epsfig}
\usepackage{amsmath,amssymb,amsthm}
\usepackage{enumerate}
\usepackage{algorithm}
\usepackage{url}
\usepackage{epstopdf}
\usepackage{graphicx}
\usepackage{pdflscape}

\DeclareMathVersion{mymath}
\SetSymbolFont{operators}{mymath}{OT1}{pplx}{b}{n}

\newenvironment{Exx}{\mathversion{mymath}}{\par\noindent}

\def\A{{\mathbf A}}
\def\PP{{\cal P}}
\def\al{\alpha}
\def\wt{\widetilde}
\def\OO{{\cal O}}

\def\proof{\medskip\noindent{\bf Proof.} }
\def\qed{\hfill $\Box$}
\newcommand{\ds}{\displaystyle}

\newtheorem{lemma}{Lemma}[section]
\newtheorem{proposition}[lemma]{Proposition}
\newtheorem{corollary}[lemma]{Corollary}
\newtheorem{theorem}[lemma]{Theorem}
\newtheorem{remark}[lemma]{Remark}
\newtheorem{conjecture}[lemma]{Conjecture}

\begin{document}

\title{Estimates for the best constant in a Markov $L_2$--inequality
       with the assistance of computer algebra}

\author{G.\,Nikolov, \ R.\,Uluchev}

\date{}
\maketitle


\begin{abstract}
We prove two-sided estimates for the best (i.e., the smallest possible)
constant $\,c_n(\alpha)\,$ in the Markov inequality
$$
  \|p_n'\|_{w_\alpha} \le c_n(\alpha) \|p_n\|_{w_\alpha}\,, \qquad p_n \in {\cal P}_n\,.
$$
Here, ${\cal P}_n$ stands for the set of algebraic polynomials of degree
$\le n$, $\,w_\alpha(x) := x^{\alpha}\,e^{-x}$, $\,\alpha > -1$, is the
Laguerre weight function, and  $\|\cdot\|_{w_\alpha}$ is the associated
$L_2$-norm,
$$
  \|f\|_{w_\alpha} = \left(\int_{0}^{\infty} |f(x)|^2 w_\alpha(x)\,dx\right)^{1/2}\,.
$$
Our approach is based on the fact that $\,c_n^{-2}(\alpha)\,$ equals
the smallest zero of a polynomial $\,Q_n$, orthogonal with respect
to a measure supported on the positive axis and defined by an
explicit three-term recurrence relation. We employ computer algebra
to evaluate the seven lowest degree coefficients of $\,Q_n\,$ and to
obtain thereby bounds for $\,c_n(\alpha)$. This work is a continuation
of a recent paper \cite{ns17a}, where estimates for $\,c_n(\alpha)\,$
were proven on the basis of the four lowest degree coefficients of $\,Q_n$.

\medskip
{\bf Keywords:} Markov type inequalities, orthogonal polynomials,
                Laguerre weight function, three-term recurrence relation,
                computer algebra.

{\bf 2000 Math.\ Subject Classification:} 41A17
\end{abstract}


\bigskip
\section{Introduction and statement of the results}

Throughout this paper $\,\PP_n\,$ will stand for the set of
algebraic polynomials of degree at most $~n$, assumed, without loss
of generality, with real coefficients. Let $\,w_\al(x) :=
x^{\al}\,e^{-x}$, where $\al
> -1$, be the Laguerre weight function, and $\,\|\cdot\|_{w_\al}\,$
be the associated $L_2$-norm,
$$
  \|f\|_{w_\al} = \left(\int_{0}^{\infty} |f(x)|^2 w_\al(x)\,dx\right)^{1/2}\,.
$$
 We study the best constant $\,c_n(\alpha)\,$
in the Markov inequality in this norm
\begin{equation} \label{e1.1}
  \|p_n'\|_{w_\al} \le c_n(\alpha) \|p_n\|_{w_\al}\,,\qquad p_n \in \PP_n\,,
\end{equation}
namely the constant
$$
  c_n(\alpha) := \sup_{p_n \in \PP_n} \frac{\|p_n'\|_{w_\al}}{\|p_n\|_{w_\al}}\,.
$$

Before formulating our results, let us give a brief account on the
results known so far.

\smallskip
It is only the case $\,\al=0\,$ where the best Markov constant is
known, namely, Tur\'{a}n \cite{pt60} proved that
$$
  c_n(0) = \Big(2\sin\frac{\pi}{4n+2}\Big)^{-1}\,.
$$
D\"{o}rfler \cite{pd91} showed that $c_n(\alpha)=\OO(n)$ for every fixed
$\,\al>-1\,$ by proving the estimates
\begin{eqnarray}
  \hspace*{-12mm} & & c_n^2(\alpha)\geq\frac{n^2}{(\al+1)(\al+3)}
      + \frac{(2\al^2+5\al+6)\,n}{3(\al+1)(\al+2)(\al+3)}
      + \frac{\al+6}{3(\al+2)(\al+3)}\,, \label{e1.2}\vspace*{2mm}\\
  \hspace*{-12mm} & & c_n^2(\alpha)\leq\frac{n(n+1)}{2(\al+1)}\,, \label{e1.3}
\end{eqnarray}
see \cite{pd02} for a more accessible source. In the same paper,
\cite{pd02}, D\"{o}rfler proved for the asymptotic constant
\begin{equation} \label{e1.4}
  c(\alpha) := \lim_{n\rightarrow\infty}\frac{c_n(\alpha)}{n}\,,
\end{equation}
that
\begin{equation} \label{e1.4Bes}
  c(\alpha) = \frac{1}{j_{(\alpha-1)/2,1}}\,,
\end{equation}
where $j_{\nu,1}$ is the first positive zero of the Bessel function
$J_{\nu}(z)$\,.

Nikolov and Shadrin obtained in \cite{ns17a} the following result:

\medskip\noindent
\textbf{Theorem A (\cite[Theorem 1]{ns17a}).}
\emph{For all $\,\alpha>-1\,$ and $\,n\in \mathbb{N}\,$, $\,n\geq 3\,$,
the best constant $\,c_n(\alpha)\,$ in the Markov inequality \eqref{e1.1}
admits the estimates
\begin{equation} \label{e1.5}
  \frac{2 \big(n+\frac{2\alpha}{3}\big)\big(n-\frac{\alpha+1}{6}\big)}
       {(\alpha+1)(\alpha+5)}
  < c_n^2(\alpha) <
  \frac{\big(n+1\big) \big(n+\frac{2(\alpha+1)}{5}\big)}
       {(\alpha+1)\big[(\alpha+3)(\alpha+5)\big]^{1/3}}\,,
\end{equation}
where for the left-hand inequality it is additionally assumed that
$\,n>(\alpha+1)/6$\,.}

\medskip
Theorem~A implies some inequalities for the asymptotic Markov
constant $\,c(\alpha)\,$ and, through \eqref{e1.4Bes}, inequalities for
$\,j_{\nu,1}\,$, the first positive zero of the Bessel function
$\,J_{\nu}\,$ (see \cite[Corollaries~1,\,3]{ns17a}). It was also
shown in \cite[Theorem~2]{ns17a} that
$\,c(\alpha)=\OO(\al^{-1})\,$, which indicates that the upper estimate
for $\,c_n(\alpha)\,$ in Theorem~A, though rather good for moderate
$\,\al\,$, is not optimal.

\smallskip
In a recent paper \cite{ns17b} Nikolov and Shadrin proved an upper
bound for $\,c_n(\alpha)\,$ which is of the correct order with respect
to both $\,n\,$ and $\,\al\,$ as they tend to infinity.

\medskip\noindent
\textbf{Theorem B (\cite[Theorem 1.1]{ns17b}).}
\emph{ For all $\,n\in \mathbb{N}\,$, $\,n\geq 3\,$, the best constant
$\,c_n(\alpha)\,$ in the Markov inequality \eqref{e1.1} satisfies
the inequality
\begin{equation}
  c_n^2(\alpha) \leq \frac{4n\big(n+2+\frac{3(\al+1)}{4}\big)}{\al^2+10\al+8}\,,
  \qquad \al\geq 2\,. \label{e1.6}
\end{equation}}

\medskip
As a consequence of Theorem~B and D\"{o}rfler's lower bound \eqref{e1.2}
for $\,c_n(\alpha)\,$ Nikolov and Shadrin showed that
$$
  c_n^2(\alpha) \asymp \frac{n(n+\al+3)}{(\al+1)(\al+8)}\,,\qquad n\geq 3,\ \al\geq 2\,.
$$

\noindent
\textbf{Corollary C (\cite[Corollary 1.1]{ns17b}).}
\emph{ For all $\,\al\geq 2\,$ and  $\,n\geq 3\,$ the best constant
$\,c_n(\alpha)\,$ in the Markov inequality \eqref{e1.1} satisfies
\begin{equation} \label{e1.7}
  \frac{2n(n+\al+3)}{3(\al+1)(\al+8)} \leq c_n^2(\alpha) \leq
  \frac{4n(n+\al+3)}{(\al+1)(\al+8)}\,.
\end{equation}}

\smallskip
In addition, Nikolov and Shadrin found the limit value of
$\,(\al+1)c_n^2(\alpha)\,$ as $\,\al\to -1$, and proved
asymptotic inequalities for $\,\al\,c_n^2(\alpha)\,$ as
$\,\al\rightarrow\infty\,$.

\medskip\noindent
\textbf{Corollary D (\cite[Corollary 1.2]{ns17b}).}
\emph{ The best constant $\,c_n(\alpha)\,$ in the Markov inequality
\eqref{e1.1} satisfies:}

\noindent
\emph{(i)} \qquad $\ds \lim_{\al\rightarrow -1} (\al+1)c_n^2(\alpha)=\frac{n(n+1)}{2}\,$;

\noindent
\emph{(ii)} \qquad $\ds \frac{2n}{3}\leq\lim_{\al\to \infty}\al\,c_n^2(\alpha)\leq 3n\,$.

\pagebreak
A combination of Theorem~A and Theorem~B implies some inequalities for the
asymptotic Markov constant \eqref{e1.4}:

\medskip\noindent
\textbf{Corollary E (\cite[Corollary 1.2]{ns17b}).}
\emph{ The asymptotic Markov constant
$\,c(\alpha) = \lim\limits_{n\to\infty} \dfrac{c_n(\alpha)}{n}\,$ satisfies the
inequalities
$$
  \frac{2}{(\al+1)(\al+5)} < c^2(\alpha) <
  \begin{cases}
    \, \dfrac{1}{(\al+1)\sqrt[3]{(\al+3)(\al+5)}}\,, & \quad -1<\al\leq \al^{*}\,, \vspace*{1mm}\\
    \, \dfrac{4}{\al^2+10\al+8}\,, &\quad \al>\al^{*}\,,
 \end{cases}
$$
where $\al^{*}\approx 43.4$\,.}

The ratio of the upper and the lower bound for $\,c(\alpha)\,$ in Corollary~E
is less than $\,\sqrt{2}\,$ for all $\,\al>-1\,$.

\smallskip
In this paper we investigate the best Markov constant $\,c_n(\alpha)\,$
following the approach from \cite{ns17a}. It is known (see
Proposition~\ref{p2.1} below) that $\,c_n^{-2}(\alpha)\,$ is equal to
the smallest zero of a polynomial $\,Q_n\,$, which is orthogonal
with respect to a measure supported on $\,\mathbb{R}_{+}\,$. Since
$\,\{Q_n\}_{n\in \mathbb{N}}\,$ are defined by an explicit
three-term recurrence relation, one can evaluate (at least
theoretically) as many coefficients of $\,Q_n\,$ as necessary. With
the assistance of Wolfram's \textsl{Mathematica} we find the seven
lowest degree coefficients of the polynomial $\,Q_n\,$, and thereby
the six highest degree coefficients of $\,R_n\,$, the monic
polynomial reciprocal to $\,Q_n\,$. Then we apply a simple technique
for estimating the largest zero $\,x_n\,$ of $\,R_n\,$ on the basis
of its $k$ highest degree coefficients, $\,3\leq k\leq 6\,$, thus
obtaining lower and upper bounds for $\,c_n^{2}(\alpha)\,$. Our main
result in this paper is:

\begin{theorem} \label{t1.1}
For $\,3\leq k\leq 6\,$ and for all $\,n\geq k\,$, the best constant
$\,c_n(\alpha)\,$ in the Markov inequality \eqref{e1.1} admits the
estimates
\begin{equation} \label{e1.8}
  \underline{c}_{\,n,k}(\alpha) \leq c_n(\alpha) \leq
  \overline{c}_{n,k}(\alpha)\,, \qquad \al>-1\,,
\end{equation}
where
\begin{eqnarray}
  & & \underline{c}_{\,n,3}^2(\alpha) =
      \frac{2\,n\big(n+\frac{3(\al+1)}{8}\big)}{(\al+1)(\al+5)}\,,
      \label{e1.9}\\
  & & \overline{c}_{n,3}^{\:2}(\alpha) =
      \frac{(n+1)\big(n+\frac{2(\al+1)}{5}\big)}
           {(\al+1)\big[(\al+3)(\al+5)\big]^{1/3}}\,, \label{e1.10}\\
  & & \underline{c}_{\,n,4}^2(\alpha) =
      \frac{(5\al+17)\,n\big(n+\frac{8(\al+1)}{25}\big)}{2(\al+1)(\al+3)(\al+7)}\,,
      \label{e1.11}\\
  & & \overline{c}_{n,4}^{\:2}(\alpha) =
      \frac{(5\al+17)^{1/4}(n+1)\big(n+\frac{3(\al+1)}{7}\big)}
           {(\al+1)(\al+3)^{1/2}\big[2(\al+5)(\al+7)\big]^{1/4}}\,,
      \label{e1.12}\\
  & & \underline{c}_{\,n,5}^2(\alpha) =
      \frac{2(7\al+31)n\big(n+\frac{25(\al+1)}{84}\big)}{(\al+1)(\al+9)(5\al+17)}\,,
      \label{e1.13}\\
  & & \overline{c}_{n,5}^{\:2}(\alpha) =
      \frac{(7\al+31)^{1/5}(n+1)\big(n+\frac{4(\al+1)}{9}\big)}
           {(\al+1)(\al+3)^{2/5}\big[(\al+5)(\al+7)(\al+9)\big]^{1/5}}\,,
      \label{e1.14}\\
  & & \underline{c}_{\,n,6}^2(\alpha) =
      \frac{\big(21\al^3+299\al^2+1391\al+2073\big)
      n\big(n+\frac{2(\al+1)}{7}\big)}{(\al+1)(\al+3)(\al+5)(\al+11)(7\al+31)}\,,
      \label{e1.15}\\
  & & \overline{c}_{n,6}^{\:2}(\alpha) =
      \frac{\big(21\al^3+299\al^2+1391\al+2073\big)^{1/6}(n+1)\big(n+\frac{5(\al+1)}{11}\big)}
           {(\al+1)(\al+3)^{1/2}(\al+5)^{1/3}\big[(\al+7)(\al+9)(\al+11)\big]^{1/6}}\,.
      \label{e1.16}
\end{eqnarray}
\end{theorem}

\pagebreak
\begin{remark} \label{r1.2}
For $\,3\leq k\leq 6$, the pair
$\,\big(\underline{c}_{\,n,k}(\alpha),
\overline{c}_{\,n,k}(\alpha)\big)\,$ of bounds for
$\,c_{n}(\alpha)\,$ is deduced with the use of the $\,k\,$ highest
degree coefficients of the polynomial $\,R_n\,$ (and \eqref{e1.10}
is the upper bound obtained in \cite{ns17a}). Generally, the bounds
for $\,c_{n}(\alpha)\,$ obtained with larger $\,k\,$ are better,
although some exceptions are observed for small $\,n\,$ and
$\,\al\,$.

Clearly, inequalities \eqref{e1.8} imply bounds for the asymptotic
Markov constant $\,c(\alpha)\,$. Here, it is not difficult
to prove that the larger $\,k\,$, the better the implied lower and
upper bounds for $\,c(\alpha)\,$, hence the best bounds for
$\,c(\alpha)\,$ are obtained from \eqref{e1.8} with $\,k=6\,$.
\end{remark}

Thus, Theorem~\ref{t1.1} yields an improvement of the estimates for
the asymptotic Markov constant $\,c(\alpha)\,$ in Corollary~E.

\begin{corollary} \label{c1.3}
The asymptotic Markov constant
$\,c(\alpha) = \lim\limits_{n\to\infty}n^{-1} c_n(\alpha)\,$ satisfies the
inequalities
$$
  \underline{c}(\alpha) < c(\alpha) < \overline{c}(\alpha)\,,
$$
where
$$
  \underline{c}^2(\alpha) := \frac{21\al^3+299\al^2+1391\al+2073}
                                  {(\al+1)(\al+3)(\al+5)(\al+11)(7\al+31)}
$$
and
$$
  \overline{c}^{\:2}(\alpha) :=
  \begin{cases}
    \dfrac{\big(21\al^3+299\al^2+1391\al+2073\big)^{1\!/6}}
             {(\al+1)(\al+3)^{1\!/2}(\al+5)^{1\!/3}\big[(\al+7)(\al+9)(\al+11)\big]^{1\!/6}}\,,
    & \!\! -1\! <\! \al\leq \al^{\star}\,,\vspace*{1mm}\\
    \dfrac{4}{\al^2+10\al+8}\,, & \!\! \al>\al^{\star},
\end{cases}
$$
with $\,\al^{\star}\approx 172\,$.
\end{corollary}

It is worth noticing that the ratio of the upper and the lower bound
for $\,c(\alpha)\,$ in Corollary~\ref{c1.3} does no exceed
$\,\frac{2\sqrt{3}}{3}\approx 1.1547\,$ for all $\,\al>-1\,$.
\smallskip

Theorem~\ref{t1.1}, in particular inequality \eqref{e1.15}, implies
an improvement of the lower bound in Corollary~D(ii).

\begin{corollary}\label{c1.4}
The best constant $\,c_n(\al)\,$ in the Markov inequality
\eqref{e1.1} satisfies:
$$
\frac{6n}{7}\leq\lim_{\al\to \infty}\al\,
       c_n^2(\alpha)\leq 3n\,.
$$
\end{corollary}

\smallskip
The rest of the paper is organized as follows. Sect.~2 contains some
preliminaries. In Sect.~2.1 we characterize the squared best Markov
constant as the largest zero of an $n$-th degree monic polynomial
$\,R_n\,$  with positive roots, and propose a recursive procedure
for the evaluation of its coefficients (Proposition~\ref{p2.2}).
Two-sided estimates for the largest zero of polynomials with only
positive roots in terms of few of their coefficients are proposed in
Sect.~2.2 (Proposition~2.3). The assisted by Wolfram's
\textsl{Mathematica} proof of our results  is given in Sect.~3. In
Sect.~4 we give some final remarks and conclusions, and formulate
two conjectures concerning the asymptotic behaviour of the best
Markov constant and the coefficients of the characteristic
polynomial $\,R_n\,$.


\bigskip
\setcounter{equation}{0}
\section{Preliminaries}

\subsection{An orthogonal polynomial related to \boldmath{$c_n(\alpha)$}}

It is well-known that the squared best constant in a Markov-type inequality
in $L_2$-norm is equal to the largest eigenvalue of a related positive
definite $\,n\times n\,$ matrix $\,\A_n\,$, thus the problem of finding
the best Markov constant is equivalent to evaluating the largest eigenvalue
of $\,\A_n$. Perhaps, a less known fact is that for a wide class of
$L_2$-norms, the inverse matrix $\,\A_n^{-1}\,$ is tri-diagonal,
see \cite[Sect. 2]{an17}. In the particular case of the
$L_2$-norm induced by the Laguerre weight function $\,w_{\al}\,$
this connection is given by the following proposition:

\begin{proposition}[\textbf{\cite[p. 85]{pd02}}] \label{p2.1}
The quantity $\,c_n^{-2}(\alpha)\,$ is equal to the smallest zero
of the polynomial $\,Q_n(x)=Q_n(x,\alpha)\,$, which is defined
recursively by
\begin{eqnarray*}
  & & Q_{n+1}(x)=(x-d_n)Q_n(x)-\lambda_n^2 Q_{n-1}(x),\quad n\ge 0\,;\\
  & & Q_{-1}(x):=0,\ \ Q_0(x):=1\,; \\
  & & d_0:=1+\alpha,\ \ d_n:=2+\frac{\alpha}{n+1}\,,\quad n\ge 1\,;\\
  & & \lambda_0>0\ \ \text{{\rm arbitrary}},\
      \lambda_n^2:=1+\frac{\alpha}{n}\,,\quad n\ge 1\,.
\end{eqnarray*}
\end{proposition}

By Favard's theorem, for any $\,\alpha>-1\,$,
$\,\{Q_n(x,\alpha)\}_{n=0}^{\infty}\,$ form a system of monic orthogonal
polynomials. Since $\,Q_n\,$ is the characteristic polynomial of the
inverse of a positive definite matrix (which is also positive
definite), it follows that all the zeros of $\,Q_n\,$ are positive
(and distinct). Consequently, $\,\{Q_n\}_{n=0}^{\infty}\,$ are
orthogonal with respect to a measure supported on $\mathbb{R}_{+}$.

\smallskip
By Proposition~\ref{p2.1}, we have
\begin{eqnarray}
  \hspace*{-10mm} & & Q_{n+1}(x) = \Big(x-2-\frac{\alpha}{n+1}\Big)Q_n(x)-\Big(1+\frac{\al}{n}\Big)Q_{n-1}(x)\,,
      \quad n\geq 1\,, \label{e2.1}\\
  \hspace*{-10mm} & & Q_0(x) = 1\,,\quad Q_1(x)=x-\al-1\,. \label{e2.2}
\end{eqnarray}

If we write $\,Q_n\,$ in the form
$$
  Q_n(x) = x^n-a_{n-1,n}\,x^{n-1}+a_{n-2,n}\,x^{n-2}-\cdots+(-1)^n\,a_{0,n}\,,
$$
then
\begin{equation} \label{e2.3}
  a_{0,n} = \binom{n+\alpha}{n}\,, \qquad n\in \mathbb{N}_0\,,
\end{equation}
with the convention that the right-hand side is equal to $1$ for $\,n=0\,$.
The proof is by induction with respect to $n$. For $\,n=0,\,1\,$, \eqref{e2.3}
follows from \eqref{e2.2}. Assuming \eqref{e2.3} is true for all
$\,m\leq n\,$, we verify it for $\,m=n+1\,$ by putting $\,x=0\,$ in
\eqref{e2.1} and using the induction hypothesis:
\begin{align*}
  (-1)^{n+1}a_{0,n+1}
      & = \Big(2+\frac{\alpha}{n+1}\Big)(-1)^{n+1}\binom{n+\alpha}{n} +
          \Big(1+\frac{\al}{n}\Big)(-1)^{n}\binom{n-1+\alpha}{n-1} \\
      & = (-1)^{n+1}\binom{n+1+\alpha}{n}\,.
\end{align*}

Now, instead of  $\,\{Q_n\}_{n=0}^{\infty}\,$, we consider the
sequence of orthogonal polynomials $\,\{\wt{Q}_n\}_{n=0}^{\infty}\,$
normalised so that $\,\wt{Q}_n(0)=1$\,, $\,n\in \mathbb{N}_0\,$,
i.e.,
$$
  Q_n(x) = (-1)^n \binom{n+\alpha}{n}\wt{Q}_n(x)\,,
  \qquad n\in\mathbb{N}_0\,.
$$
It follows from \eqref{e2.1} and \eqref{e2.2} that
$\,\{\wt{Q}_n\}_{n\in \mathbb{N}_0}\,$ are determined by
\begin{eqnarray}
  \hspace*{-10mm} & & \Big(1+\frac{\al}{n+1}\Big)\wt{Q}_{n+1}(x) =
      \Big(2+\frac{\al}{n+1}-x\Big)\wt{Q}_n(x)-\wt{Q}_{n-1}(x)\,,
      \quad n\geq 1\,, \label{e2.4}\\
  \hspace*{-10mm} & & \wt{Q}_0(x) = 1\,,\quad \wt{Q}_1(x)=1-\frac{x}{\al+1}\,. \label{e2.5}
\end{eqnarray}

Writing $\,\wt{Q}_n\,$ in the form
$$
  \wt{Q}_n(x) = 1-A_{1,n}\,x+A_{2,n}\,x^2-\cdots+(-1)^nA_{n,n}\,x^n
$$
and rewriting \eqref{e2.4} as
$$
  \wt{Q}_{n+1}(x)-\wt{Q}_n(x)=\frac{n+1}{n+\al+1}\big(\wt{Q}_{n}(x)
  -\wt{Q}_{n-1}(x)\big)+\frac{n+1}{n+\al+1}\,x\,\wt{Q}_{n}(x)\,,
  \quad n\in \mathbb{N}\,,
$$
we deduce the following recurrence relation for the evaluation of the
coefficients $\,\{A_{i,m}\}\,$:
\begin{equation}
  \begin{aligned} \label{e2.6}
  \hspace*{-2mm} & A_{i,n+1}-A_{i,n}=\frac{n+1}{n+\al+1}\big(A_{i,n}-A_{i,n-1}\big)+
      \frac{n+1}{n+\al+1}\,A_{i-1,n}\,,\ \ n\geq k\geq 1\,, \\
  \hspace*{-2mm} & \text{with } A_{0,n}=1\ \text{ and }\ A_{1,1}=\frac{1}{\al+1}\,.
 \end{aligned}
\end{equation}

Since, by Proposition~\ref{p2.1}, $\,c_n^{-2}(\alpha)\,$ is equal
to the smallest zero of $\,\wt{Q}_n\,$, it follows that
$\,c_n^2(\alpha)\,$ equals the largest zero of the reciprocal
polynomial of $\,\wt{Q}_n\,$,
\begin{equation} \label{e2.6rn}
  R_n(x)=x^{n}\,\wt{Q}_n(1/x)\,.
\end{equation}

The above observations allow us to reformulate
Proposition~\ref{p2.1} in the following equivalent form:

\begin{proposition} \label{p2.2}
The squared best Markov constant $\,c_n^2(\alpha)\,$ is equal to
the largest zero of the polynomial
\begin{equation} \label{e2.7}
  R_n(x) = x^n-A_{1,n}\,x^{n-1}+A_{2,n}\,x^{n-2}-\cdots+(-1)^n A_{n,n}\,.
\end{equation}

The coefficients of $\,R_n\,$ are evaluated recursively by the
following procedure:
\begin{itemize}
  \item $\,A_{1,1}=\frac{1}{\al+1}$\,;
  \item Set $\,A_{0,m}=1\,$, $\,m=0,\ldots,n\,$;
  \item For $i=1$ to $n$:
        \begin{enumerate}
          \item Find the sequence $\,\{D_{i,m}\}_{m=i-1}^n\,$
                as solution of the recurrence equation
                \begin{equation} \label{e2.8}
                  D_{i,m+1}=\frac{m+1}{m+\al+1}\,D_{i,m}+\frac{m+1}{m+\al+1}\,A_{i-1,m}
                \end{equation}
                with the initial condition $\,D_{i,i-1}=0\,$;
          \item Evaluate
                \begin{equation} \label{e2.8akn}
                  A_{i,n}=\sum_{m=i}^n D_{i,m}\,.
                \end{equation}
        \end{enumerate}
\end{itemize}
\end{proposition}

\subsection{Polynomials with positive roots: bounds for the largest zero}

Let $\,P\,$ be a monic polynomial of degree $\,n\,$ with zeros
$\,\{x_i\}_{i=1}^{n}\,$,
$$
  P(x) = \prod_{i=1}^{n}(x-x_i)=x^n-b_1\,x^{n-1}+b_2\,x^{n-2}-\cdots+(-1)^n b_n\,.
$$

The coefficients $\,b_r=b_r(P)\,$, $\,r=1,\ldots,n\,$, are given by
the elementary symmetric functions of $\,\{x_i\}_{i=1}^n\,$,
$$
  b_r = s_r = s_r(P) =
  \sum_{1\leq i_1<i_2<\cdots<i_r\leq n} x_{i_1}x_{i_2}\cdots x_{i_r}\,,
  \qquad r=1,\ldots,n\,.
$$
It is well known that the elementary symmetric functions
$\,\{s_r\}\,$ and the Newton functions (sums of powers of $\,x_i\,$)
$$
  p_r = p_r(P) = \sum_{i=1}^n x_i^r\,, \qquad r=1,2,3,\ldots\,,
$$
are connected by the Newton identities:
\begin{alignat}{3}
   & p_r+\sum_{i=1}^{r-1}(-1)^{i}p_{r-i}\,s_i+(-1)^r r\,s_r=0\,, \qquad
       && \text{if } \ 1\leq r\leq n\,, \label{e2.9}\\
   & p_r+\sum_{i=1}^{n}(-1)^{i}p_{r-i}\,s_i=0\,, \qquad
       && \text{if } \ r>n\,. \label{e2.10}
\end{alignat}
For a proof, see e.g. \cite{vw49} or \cite{dm92}.

\smallskip
Our interest in the Newton functions is motivated by the fact that they
provide tight bounds for the largest zero of a polynomial whose roots are
all positive. For any such polynomial $\,P\,$, we set
$$
  \ell_k(P) := \frac{p_k(P)}{p_{k-1}(P)}\,, \qquad
  u_k(P) := \big[p_k(P)\big]^{1/k}\,, \qquad k\in \mathbb{N}\,,
$$
with the convention that $\,p_0(P):=\deg(P)\,$.

\begin{proposition} \label{p2.3}
Let
$\,P(x)=x^n-b_{1}\,x^{n-1}+b_{2}\,x^{n-2}-\cdots+(-1)^{n-1}b_{n-1}\,x+(-1)^{n}b_n\,$
be a polynomial with positive zeros $\,x_1\le x_2\le\cdots\le x_n\,$.

Then the largest zero $\,x_n\,$ of $\,P\,$  satisfies the
inequalities
\begin{equation} \label{e2.11}
  \ell_k(P)\leq x_n < u_k(P)\,, \qquad k\in \mathbb{N}\,.
\end{equation}
Moreover, the sequence $\,\{\ell_k(P)\}_{k=1}^{\infty}\,$ is
monotonically increasing, the sequence
$\,\{u_k(P)\}_{k=1}^{\infty}\,$ is monotonically decreasing, and
\begin{equation} \label{e2.12}
  \lim_{k\to\infty} \ell_k(P) = \lim_{k\to\infty} u_k(P) = x_n\,.
\end{equation}
\end{proposition}

\proof For $\,i=1,\ldots,n-1\,$, we set $\,a_i:=\frac{x_i}{x_n}\,$, then
$\,0<a_i\leq 1\,$. Now both inequalities \eqref{e2.11} and the limit
relations \eqref{e2.12} readily follow from the representations
$$
  \ell_k(P) = \frac{a_1^k+\cdots+a_{n-1}^k+1}{a_1^{k-1}+\cdots+a_{n-1}^{k-1}+1}\,x_n\,,
  \qquad
  u_k(P) = \big(a_1^k+\cdots+a_{n-1}^k+1\big)^{1/k} x_n\,.
$$

The monotonicity of the sequence $\,\{\ell_k(P)\}_{k=1}^{\infty}\,$
follows easily from Cauchy-Bouniakowsky's inequality. Indeed, we
have
$$
  \Big(\sum_{i=1}^{n}x_i^{k}\Big)^2 =
  \Big(\sum_{i=1}^{n}x_i^{\frac{k-1}{2}}x_i^{\frac{k+1}{2}}\Big)^2
  \leq
  \Big(\sum_{i=1}^{n}x_i^{k-1}\Big)\Big(\sum_{i=1}^{n}x_i^{k+1}\Big)\,,
$$
whence $\,p_k^2(P)\leq p_{k-1}(P)\,p_{k+1}(P)\,$, and consequently
$$
  \ell_k(P) = \frac{p_k(P)}{p_{k-1}(P)} \leq
  \frac{p_{k+1}(P)}{p_k(P)} = \ell_{k+1}(P)\,.
$$

To prove monotonicity of the sequence
$\,\{u_k(P)\}_{k=1}^{\infty}\,$, we recall that $\,0< a_i\leq 1\,$
and therefore $\,a_i^{k+1}\leq a_i^k\,$. We have
$$
  \big(a_1^{k+1}+\cdots+a_{n-1}^{k+1}+1\big)^{1/(k+1)} \!<
  \big(a_1^{k+1}+\cdots+a_{n-1}^{k+1}+1\big)^{1/k} \!\leq
  \big(a_1^k+\cdots+a_{n-1}^k+1\big)^{1/k},
$$
which yields
$$
  u_{k+1}(P) < u_k(P)\,.
$$
\qed


\bigskip
\setcounter{equation}{0}
\section{Computer algebra assisted proof of the results}

Here we give the algorithms, the source code and the results of the
computer algebra assisted proof of estimates \eqref{e1.9}-\eqref{e1.16} in
Theorem~\ref{t1.1}. While the case $\,k=3\,$ and to a certain extent
$\,k=4\,$ could be studied by hand, it seems impossible to provide
similar calculations for larger~$\,k\,$. We implement the idea from
\cite{ns17a} for estimating $\,c_n(\alpha)\,$ using $\,k=3\,$ highest
degree coefficients of the polynomial $\,R_n(x)\,$ and with the assistance
of Wolfram's \textsl{Mathematica} v. 10 software we investigate the cases
$\,k=4,5,6\,$, as well. Software based on the algorithms described below
failed with calculations for $\,k>6\,$.

For simplicity sake, henceforth we write the polynomial $\,R_n\,$
from \eqref{e2.6rn} and \eqref{e2.7} in the form
$$
  R_n(x) = x^n - b_1 x^{n-1} + b_2 x^{n-2} + \cdots + (-1)^nb_n\,.
$$


\subsection{Lower bounds for \boldmath{$c_n(\alpha)$}} \label{sect3.1}

We apply Proposition~\ref{p2.3} to estimate the largest zero
$\,x_n=c_n^2(\al)\,$ of the polynomial $\,R_n(x)\,$ from below,
$$
x_n\geq \ell_k(R_n)=\frac{p_k(R_n)}{p_{k-1}(R_n)}\,,\qquad
k=3,\,4,\,5,\,6\,,
$$
and then with the help of computer algebra obtain a further
estimation of the form
$$
\ell_k(R_n)\geq c\,n(n+\sigma(\alpha+1)),
$$
with the optimal (i.e., the largest possible) constants $\,c=c(k)\,$
and $\,\sigma=\sigma(k)$.

\begin{algorithm}
\caption{\quad Estimating $\,c_n(\alpha)\,$ from below} \label{alg1}
\begin{tabular}{ll}
  {\sl Input:}  & $\,k\in\{3,4,5,6\}\,$ -- the number of the highest degree coefficients of $\,R_n(x)\,$ \\
  {\sl Step 1.} & Express the power sums $\,p_{k-1}(R_n)\,$ and $\,p_k(R_n)\,$ in terms of $\,\{b_i\}_{i=1}^k\,$ \\
  {\sl Step 2.} & Find coefficients $\,\{b_i\}_{i=1}^k\,$ in terms of $\,n\,$ and $\,\alpha\,$ using Proposition~\ref{p2.2} \\
  {\sl Step 3.} & Find a proper value $\,\sigma\,$ for parameter $\,s\,$ in $\,p_k-c\,n(n+s(\alpha+1))p_{k-1}\,$, \\
                & where $\,c\,$ is the coefficient of $\,n^2\,$ in the quotient $\,p_k/p_{k-1}\,$ \\
  {\sl Step 4.} & Represent the numerator of $\,f=p_k-c\,n(n+\sigma(\alpha+1))p_{k-1}\,$ in powers \\
                & of $\,n\,$ and $\,(\alpha+1)\,$ \\
  {\sl Step 5.} & Estimate from below the expression $\,f\,$ to prove that $\,f\ge 0\,$
\end{tabular}
\end{algorithm}

\smallskip
\underline{\sl Step 1:} \ Let $\,\{x_i\}_{i=1}^n\,$ be all the zeros
of the polynomial $\,R_n(x)\,$ from \eqref{e2.6rn}. In order to
express a power sum $\,p_r=\sum_{i=1}^r x_i^r\,$, $\,1\le r\le n$,
by $\{b_i\}_{i=1}^r\,$, we apply the direct formula
$$
  p_r=\begin{vmatrix}
         b_1 &   1 &  0 & \cdots & 0 \\
        2b_2 & b_1 &  1 & \cdots & 0 \\
        3b_3 & b_2 & b_1 & \cdots & 0 \\
        \hdotsfor{5} \\
        r b_r & b_{r-1} & b_{r-2} & \cdots & b_1
      \end{vmatrix}
$$
which easily follows from the Newton identities \eqref{e2.9}.

Below is the code of the programme and the results for
$\,k=1,\ldots,6\,$:
\begin{center}
  \epsfig{file=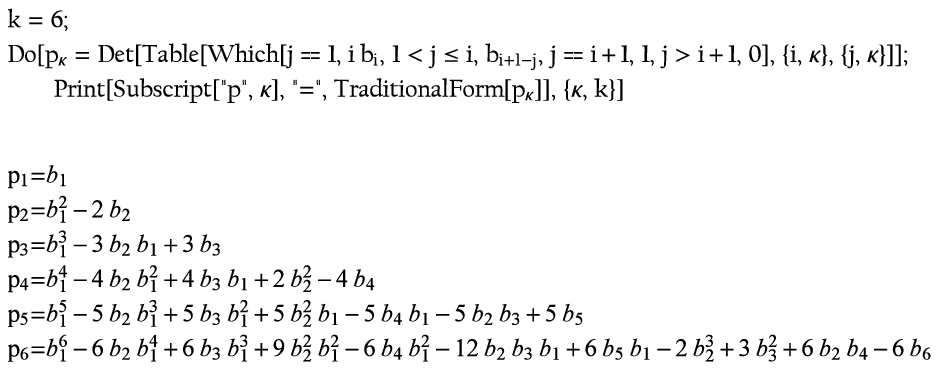}
\end{center}

\smallskip
\underline{\sl Step 2:} \
We find coefficients $\,\{b_i\}_{i=1}^k\,$ of the polynomial
$\,R_n(x)\,$ using Proposition~\ref{p2.2}. For a fixed $\,i\,$ we
firstly find a sequence solving recurrence equation \eqref{e2.8} and
then evaluate $\,b_i\,$ by \eqref{e2.8akn}.

The source and the results for $\,k=1,\ldots,6\,$ follow below:
\begin{center}
  \epsfig{file=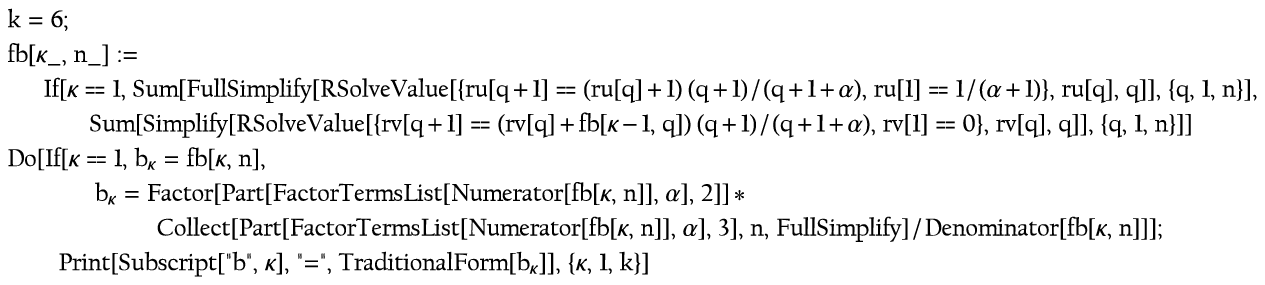}
\end{center}

\begin{center}
  \epsfig{file=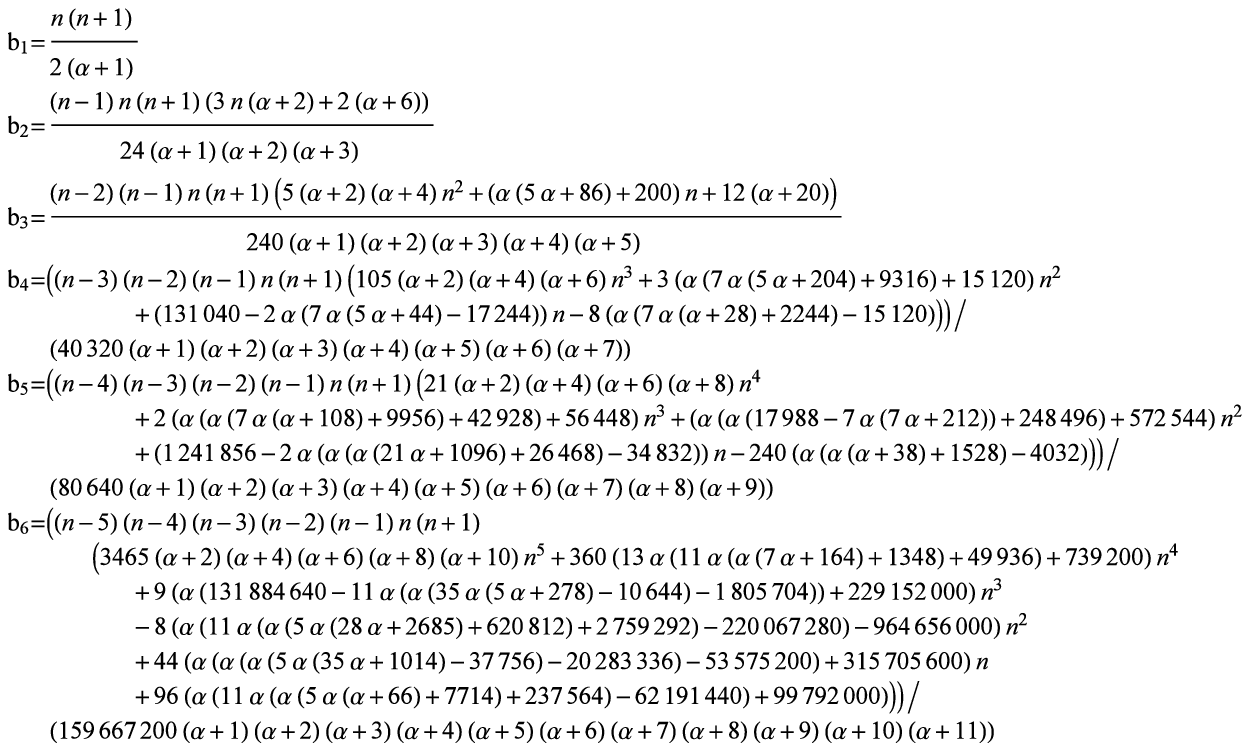}
\end{center}

\smallskip
\underline{\sl Step 3:} \
The quotient $\,p_k/p_{k-1}\,$ is a quadratic polynomial in $\,n\,$,
and we denote by $\,c\,$ its leading coefficient.

The goal of this step is to find a proper value (say $\sigma$) for
parameter $\,s\,$ in the expression
$$
  f_s=p_k-c\,n(n+s(\alpha+1))p_{k-1}\,,
$$
such that $\,f_{\sigma}\geq 0\,$ for all admissible $\,\alpha\,$ and
$\,n\,$. For a fixed $\,k\,$ quantity $\,f_s\,$ depends on
$\,\alpha\,$, $\,n\,$ and $\,s\,$. It is a polynomial of degree
$\,2k-1\,$ in $\,n\,$ and a rational function in $\,\alpha\,$. Let
us write the numerator of $\,f_s\,$ in the form
$$
  \sum_{i=1}^{2k-1} \sum_{j=0}^d \mu_{i,j}(s)(\alpha+1)^{d-j} n^{2k-i}\,.
$$
The highest order coefficients in $\,\sum_j \mu_{i,j}(s)(\alpha+1)^{d-j}\,$ are
linear functions in $\,s\,$ of the form $\,A_i-B_i s\,$, with $\,A_i>0\,$ and
$\,B_i>0\,$. We denote their zeros by $\,s_i\,$ for each $\,i\,$ and set
$\,\sigma=\min_i s_i\,$. Since we seek estimates valid for all $\,\alpha>-1\,$,
our choice of $\,\sigma\,$ guarantee that for $\,\alpha\,$ sufficiently large
the inequality $\,\sum_j \mu_{i,j}(s)(\alpha+1)^{d-j}>0\,$ holds true.

The code is as follows:
\begin{center}
  \epsfig{file=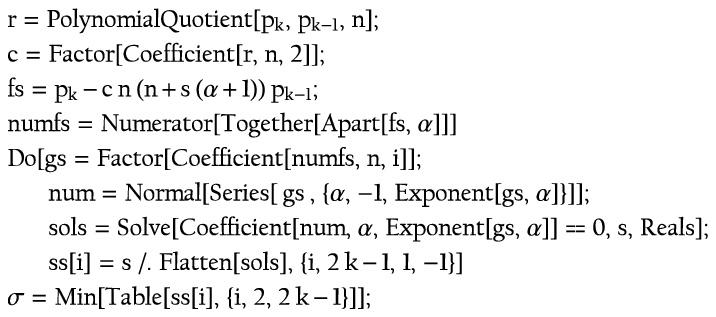}
\end{center}

Table~\ref{tab1} gives results for the optimal values of $\,c\,$ and
\,$\sigma\,$ for $\,k=3,4,5,6\,$.
\vspace{-3ex}
\begin{table}[H]
  \centering
  \caption{The optimal values of $\,c\,$ and $\,\sigma\,$ in the lower bounds for $c_n^2(\al)$.} \label{tab1}
  \medskip
  \begin{tabular}{ccc} \hline
      $k$ \qquad & $c$ & \quad $\sigma$ \quad \\[2pt] \hline
       3  \qquad & $\dfrac{2}{(\alpha+1)(\alpha+5)}$ & \quad \mbox{\rule[-4mm]{0mm}{10mm}{$\dfrac38$}} \\[8pt]
       4  \qquad & $\dfrac{5\alpha+17}{2(\alpha+1)(\alpha+3)(\alpha+7)}$ & \quad $\dfrac{8}{25}$ \\[8pt]
       5  \qquad & $\dfrac{2(7\alpha+31)}{(\alpha+1)(\alpha+9)(5\alpha+17)}$ & \quad $\dfrac{25}{84}$ \\[8pt]
       6  \qquad & \quad $\dfrac{21\alpha^3+299\alpha^2+1391\alpha+2073}{(\alpha+1)(\alpha+3)(\alpha+5)(\alpha+11)(7\alpha+31)}$
                   & \quad $\dfrac27$ \\[8pt] \hline
  \end{tabular}
\end{table}

\underline{\sl Step 4:}  We set
$$
  f=p_k-c\,n(n+\sigma(\alpha+1))p_{k-1}=:\frac{\varphi(n,\alpha)}{\psi(\alpha)}
$$
with $\,c\,$ and $\,\sigma\,$ determined in Step~3. Here,
$\,\varphi(n,\alpha)\,$ is a bivariate polynomial in $\,n\,$ and
$\,\alpha\,$, and $\,\psi(\alpha)\,$ is a polynomial in
$\,\alpha\,$. More precisely, $\,\varphi(n,\alpha)\,$ has degree
$\,2k-1\,$ in $\,n\,$, and degree $\,d\,$ in $\,\alpha\,$ which our
programme calculates for each fixed $\,k\,$.

Note that $\,\psi(\alpha)>0\,$ for $\,\alpha>-1\,$ since it is a product
of powers of $\,\alpha+j\,$, $\,j\ge 1\,$ and multipliers $\,A\alpha+B\,$,
$\,0<A<B\,$. Therefore, $\,{\rm sign}\, f = {\rm sign}\,\varphi\,$.

We expand $\,\varphi(n,\alpha)\,$ in the form
$$
  \varphi(n,\alpha) = \sum_{i=1}^{2k-1} \sum_{j=0}^d \mu_{i,j}(\alpha+1)^{d-j} n^{2k-i} =
    \begin{pmatrix} n^{2k-1} \\ n^{2k-2} \\ \vdots \\ n \end{pmatrix}^\top
    \mathbf{M}
    \begin{pmatrix} (\alpha+1)^d \\ (\alpha+1)^{d-1} \\ \vdots \\ 1 \end{pmatrix}\,,
$$
where $\,\mathbf{M}=\big(\mu_{i,j}\big)_{i=1,j=0}^{2k-1,d}\,$ and all entries
$\,\mu_{i,j}\,$ are integer numbers.

The source for computation of the matrix $\,\mathbf{M}\,$ is listed below.
\begin{center}
  \epsfig{file=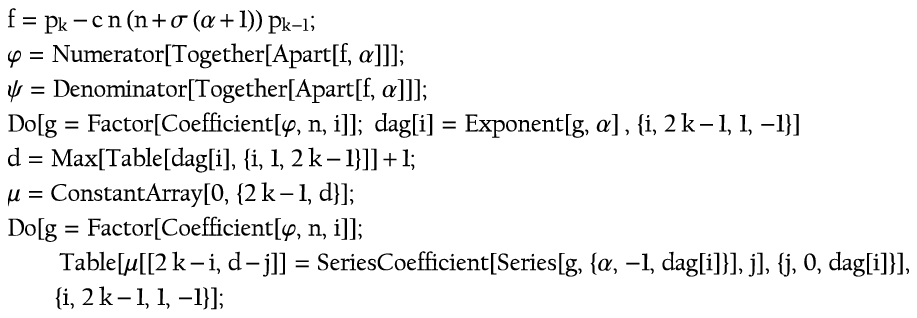}
\end{center}

If $\,\mu_{i,j}\ge 0\,$ for all $\,i,j\,$, then $\,\varphi(n,\alpha)\ge 0\,$
and $\,f\ge 0\,$ for all $\,\alpha>-1\,$ and $\,n\ge k\,$.
In a case some of coefficients $\,\mu_{i,j}<0\,$ we apply the next step of
the algorithm.

The results for $\,k=3,4,5,6\,$ are given together with the
estimates from Step~5.

\smallskip
\underline{\sl Step 5:} \ If there are coefficients
$\,\mu_{i,j}<0\,$ we need additional arguments to verify that
$\,f\ge 0\,$ for all $\,\alpha>-1\,$ and $\,n\ge k\,$. We bring into
use a new $\,(2k-1)\times (d+1)\,$ matrix $\,\mathbf{\Lambda}\,$
which elements we put initially $\,\lambda_{i,j}:=\mu_{i,j}\,$, for
$\,i=1,\ldots,2k-1\,$ and $\,j=0,\ldots,d\,$.

The procedure described below checks recursively all coefficients
$\,\lambda_{i,j}\,$ and makes the corre\-sponding estimations.
We need not introduce a new matrix after each iteration, but only
replace a pair of elements in a column of $\,\mathbf{\Lambda}\,$ with new
entries in such a manner that the value of the function
$$
  \Phi(\mathbf{\Lambda}) = \sum_{i=1}^{2k-1} \sum_{j=0}^d \lambda_{i,j}(\alpha+1)^{d-j} n^{2k-i}=
    \begin{pmatrix} n^{2k-1} \\ n^{2k-2} \\ \vdots \\ n \end{pmatrix}^\top
    \mathbf{\Lambda}
    \begin{pmatrix} (\alpha+1)^d \\ (\alpha+1)^{d-1} \\ \vdots \\ 1 \end{pmatrix}
$$
decreases. At the end of the procedure we get a matrix $\,\mathbf{\Lambda}\,$
satisfying $\,\mathbf{0}\le\mathbf{\Lambda}\le\mathbf{M}\,$ (in the sense that
$\,0\le \lambda_{i,j}\le \mu_{i,j}\,$ for all $\,i,j\,$) and therefore
$$
  0\le\Phi(\mathbf{\Lambda})\le\Phi(\mathbf{M})=\varphi(n,\alpha)\,.
$$

Suppose that $\,\lambda_{i,j}<0\,$ for some pair of indices $\,i,j\,$.
Then we set
$$
  h:=\min\{i-\eta: \lambda_{\eta,j}>0, \ 1\le \eta\le i-1\}
  \quad\text{and}\quad
  \delta:=\frac{\lambda_{i,j}}{k^{i-h}} \quad (\delta<0)\,.
$$

If $\,\lambda_{h,j}+\delta\ge 0\,$, for $\,n\ge k\,$ we have
\begin{align*}
  (\lambda_{h,j}+\delta)n^{2k-h}+0\,n^{2k-i}
   & = \Big(\lambda_{h,j}+\frac{\lambda_{i,j}}{k^{i-h}}\Big)n^{2k-h}
     = \lambda_{h,j}n^{2k-h}+\lambda_{i,j}\frac{n^{2k-h}}{k^{i-h}} \\
   & \le \lambda_{h,j}n^{2k-h}+\lambda_{i,j}\frac{n^{2k-h}}{n^{i-h}}
     = \lambda_{h,j}n^{2k-h}+\lambda_{i,j}{n^{2k-i}}\,.
\end{align*}

Otherwise, if $\,\lambda_{h,j}+\delta < 0\,$, for $\,n\ge k\,$ we have
\begin{align*}
  0\,n^{2k-h}+\big(\lambda_{h,j}k^{i-h}+\lambda_{i,j}\big)n^{2k-i}
   & = \lambda_{h,j}n^{2k-i}k^{i-h}+\lambda_{i,j}n^{2k-i}  \\
   & \le \lambda_{h,j}n^{2k-i}n^{i-h}+\lambda_{i,j}n^{2k-i} \\
   & \le \lambda_{h,j}n^{2k-h}+\lambda_{i,j}n^{2k-i}\,.
\end{align*}

So, replacing only two elements in $\,\mathbf{\Lambda}\,$,
$$
  \begin{cases}
    \lambda_{h,j}:=\lambda_{h,j}+\lfloor \delta \rfloor \ \ \qquad\, \text{and} \ \ \lambda_{i,j}:=0\,,
      & \text{ if } \lambda_{h,j}+\delta\ge 0, \\
    \lambda_{i,j}:=\lambda_{h,j}\,k^{i-h}+\lambda_{i,j} \ \ \text{and} \ \ \lambda_{h,j}:=0\,,
      & \text{ otherwise },
  \end{cases}
$$
we obtain that
$$
  \lambda_{h,j}(\alpha+1)^{d+1-j}n^{2k-h}+\lambda_{i,j}(\alpha+1)^{d+1-j}n^{2k-i}
$$
decreases for the new values of $\,\lambda_{h,j}\,$ and
$\,\lambda_{i,j}\,$, and hence $\,\Phi(\mathbf{\Lambda})\,$ also
decreases.

Applying recursively the above iteration process for
$\,i=2k-1,2k-2,\ldots,1\,$ and $\,j=0,1,\ldots,d\,$ we finally obtain a
matrix $\,\mathbf{\Lambda}\,$ satisfying
$\,\mathbf{0}\le\mathbf{\Lambda}\le\mathbf{M}\,$. Then
$\,\varphi(n,\alpha)\ge 0\,$, $\,f\ge 0\,$ and therefore
$$
c_n^2(\alpha)\ge \frac{p_k}{p_{k-1}} \ge c\,n(n+\sigma(\alpha+1))
$$
for the optimal $\,c\,$ and $\,\sigma\,$ evaluated in Step~3. For
$\,k=3,4,5,6\,$ we obtain estimates \eqref{e1.9}, \eqref{e1.11},
\eqref{e1.13}, and \eqref{e1.15}, respectively.

The following source implements the procedure described in Step~5.

\begin{center}
  \epsfig{file=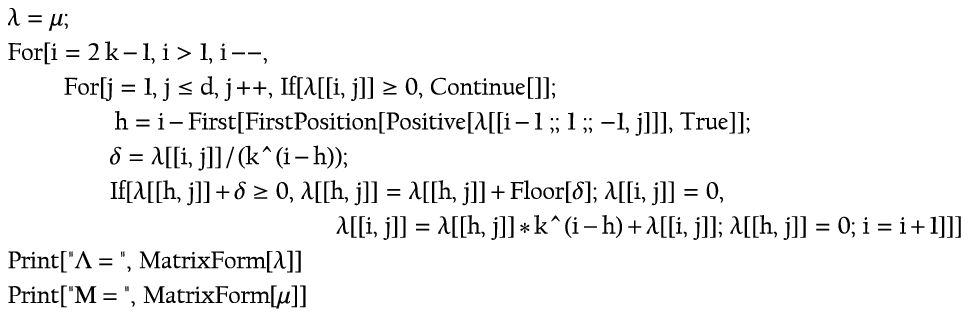}
\end{center}

Next, we give matrices $\,\mathbf{M}\,$ from Step~4 and
$\,\mathbf{\Lambda}\,$ from Step~5 obtained with {\sl Mathematica}.

\medskip
Case $\,k=3\,$:

This partial case needs a special attention as we have to assume strict
inequality $\,n>k\,$, i.e., $\,n\geq 4\,$, to obtain estimate
\eqref{e1.9}. This causes a minor modification in Step~5 of
Algorithm~\ref{alg1}, namely, replacement of $\,k^{i-h}\,$ with
$\,(k+1)^{i-h}\,$. Namely, we determine
$\,\delta:=\lambda_{i,j}/(k+1)^{i-h}\,$ and set
$$
  \begin{cases}
    \lambda_{h,j}:=\lambda_{h,j}+\lfloor \delta \rfloor \ \ \ \ \qquad\qquad\, \text{and} \ \ \lambda_{i,j}:=0\,,
      & \text{ if } \lambda_{h,j}+\delta\ge 0, \\
    \lambda_{i,j}:=\lambda_{h,j}\,(k+1)^{i-h}+\lambda_{i,j} \ \ \text{and} \ \ \lambda_{h,j}:=0\,,
      & \text{ otherwise }.
  \end{cases}
$$

\pagebreak
Matrices $\,\mathbf{M}\,$ and $\,\mathbf{\Lambda}\,$ in this case are
\begin{center}
  \begin{Exx}
    {\scriptsize
    $ \mathbf{\Lambda} =
      \left( \arraycolsep=3pt
      \begin{array}{ccccc}
        0 & 4 & -4 & 225 & 360 \\
        0 & 0 & 390 & 510 & 720 \\
       15 & 155 & 205 & 1185 & 360 \\
       15 & 270 & 495 & 900 & 0 \\
        0 & 36 & 684 & 0 & 0 \\
      \end{array}
      \right) $
    \qquad\qquad
    $ \text{M} =
      \left( \arraycolsep=3pt
      \begin{array}{ccccc}
        0 & 19 & -4 & 225 & 360 \\
        0 & -60 & 390 & 510 & 720 \\
       15 & 155 & 205 & 1185 & 360 \\
       15 & 270 & 495 & 900 & 0 \\
        0 & 36 & 684 & 0 & 0 \\
      \end{array}
      \right) $}.
  \end{Exx}
\end{center}

\noindent
Although there is a negative element of $\,\mathbf{\Lambda}\,$, from
$\,4(\alpha+1)^2-4(\alpha+1)+225\geq 0\,$ for all
$\,\alpha>-1\,$ we conclude that $\,
4(\alpha+1)^3-4(\alpha+1)^2+225(\alpha+1)+360>0\,$ and consequently
$\,\Phi(\mathbf{\Lambda})\geq 0\,$ for $\,n\geq 4$.

By a direct verification one can see that inequality \eqref{e1.9}
holds also in the case $n=k=3$.

\medskip
Case $\,k=4\,$:

\medskip
\begin{Exx}
  {\scriptsize
  $ \mathbf{\Lambda} =
    \left( \arraycolsep=3pt
    \begin{array}{ccccccccc}
      0 & 0 & 10200 & 72480 & 323700 & 1413060 & 3602340 & 4340700 & 1890000 \\
      0 & 4882 & 30891 & 359695 & 2625259 & 7966210 & 13275570 & 12707100 & 5670000 \\
      0 & 0 & 229110 & 1642830 & 6282570 & 16699200 & 24837120 & 18692100 & 5670000 \\
      2100 & 46515 & 120645 & 2404465 & 10159765 & 20026720 & 25810890 & 16625700 & 1890000 \\
      2756 & 106120 & 876330 & 2582090 & 7616630 & 17567550 & 18060000 & 6300000 & 0 \\
      0 & 11060 & 662604 & 2653840 & 6215776 & 11121880 & 7413000 & 0 & 0 \\
      0 & 0 & 0 & 1120600 & 4777900 & 3435000 & 0 & 0 & 0 \\
    \end{array}
    \right) $

  \smallskip
  $ \text{M} =
    \left( \arraycolsep=3pt
    \begin{array}{ccccccccc}
      0 & 0 & 10200 & 72480 & 323700 & 1413060 & 3602340 & 4340700 & 1890000 \\
      0 & 8715 & 30891 & 359695 & 2625259 & 7966210 & 13275570 & 12707100 & 5670000 \\
      0 & -15330 & 229110 & 1642830 & 6282570 & 16699200 & 24837120 & 18692100 & 5670000 \\
      2100 & 46515 & 120645 & 2404465 & 10159765 & 20026720 & 25810890 & 16625700 & 1890000 \\
      2800 & 106120 & 876330 & 2582090 & 7616630 & 17567550 & 18060000 & 6300000 & 0 \\
      0 & 15960 & 722904 & 2653840 & 6215776 & 11121880 & 7413000 & 0 & 0 \\
      -700 & -19600 & -241200 & 1120600 & 4777900 & 3435000 & 0 & 0 & 0 \\
    \end{array}
    \right) $ }
\end{Exx}

\medskip
Case $\,k=5\,$:

\medskip
\scalebox{0.6}{ \begin{Exx}
                  $ \mathbf{\Lambda} =
                    \left( \arraycolsep=3pt
                    \begin{array}{cccccccccccc}
                      0 & 0 & 0 & 64925 & 1064665 & 8138830 & 43256150 & 172898565 & 474925185 & 805850640 & 734423760 & 266716800 \\
                      0 & 0 & 91665 & 1204470 & 9699090 & 71280390 & 373661895 & 1241223900 & 2610599670 & 3473555400 & 2804336640 & 1066867200 \\
                      0 & 19824 & 130578 & 3408188 & 48487642 & 313463920 & 1271550350 & 3522779568 & 6544523790 & 7686433440 & 5117787360 & 1600300800 \\
                      0 & 0 & 1451982 & 16288020 & 114900450 & 672910770 & 2546690160 & 6152610870 & 9859721760 & 10218685680 & 5871579840 & 1066867200 \\
                      3675 & 128835 & 0 & 24490445 & 226233910 & 991504675 & 3153540110 & 7169071245 & 10438959825 & 9013742640 & 3935025360 & 266716800 \\
                      6027 & 381850 & 6416795 & 22404550 & 169885205 & 1005110890 & 2985302145 & 5744010510 & 7716554370 & 5584488840 & 1111320000 & 0 \\
                      0 & 52297 & 5062484 & 58263912 & 213196158 & 589342950 & 1804792500 & 3787471002 & 4038237000 & 1770703200 & 0 & 0 \\
                      0 & 0 & 0 & 15084950 & 144208510 & 409403975 & 1057769610 & 1931913900 & 1309770000 & 0 & 0 & 0 \\
                      0 & 0 & 0 & 0 & 0 & 256255650 & 690284700 & 417538800 & 0 & 0 & 0 & 0 \\
                    \end{array}
                    \right) $
                \end{Exx} }

\medskip
\scalebox{0.6}{ \begin{Exx}
                  $ \text{M} =
                    \left( \arraycolsep=3pt
                    \begin{array}{cccccccccccc}
                      0 & 0 & 0 & 64925 & 1064665 & 8138830 & 43256150 & 172898565 & 474925185 & 805850640 & 734423760 & 266716800 \\
                      0 & 0 & 91665 & 1204470 & 9699090 & 71280390 & 373661895 & 1241223900 & 2610599670 & 3473555400 & 2804336640 & 1066867200 \\
                      0 & 27804 & 130578 & 3408188 & 48487642 & 313463920 & 1271550350 & 3522779568 & 6544523790 & 7686433440 & 5117787360 & 1600300800 \\
                      0 & -39900 & 1500030 & 16288020 & 114900450 & 672910770 & 2546690160 & 6152610870 & 9859721760 & 10218685680 & 5871579840 & 1066867200 \\
                      3675 & 128835 & -240240 & 24490445 & 226233910 & 991504675 & 3153540110 & 7169071245 & 10438959825 & 9013742640 & 3935025360 & 266716800 \\
                      6125 & 381850 & 6416795 & 22404550 & 169885205 & 1005110890 & 2985302145 & 5744010510 & 7716554370 & 5584488840 & 1111320000 & 0 \\
                      0 & 77616 & 5699022 & 58263912 & 213196158 & 589342950 & 1804792500 & 3787471002 & 4038237000 & 1770703200 & 0 & 0 \\
                      -2450 & -123445 & -3055430 & 20292530 & 152590030 & 409403975 & 1057769610 & 1931913900 & 1309770000 & 0 & 0 & 0 \\
                      0 & -15750 & -636300 & -26037900 & -41907600 & 256255650 & 690284700 & 417538800 & 0 & 0 & 0 & 0 \\
                    \end{array}
                    \right) $
                \end{Exx} }

\bigskip
Case $\,k=6\,$:

\medskip\noindent
\resizebox{\linewidth}{!}{
  \begin{Exx}
    {\large $ \mathbf{\Lambda^T} = $}
    $ \left( \arraycolsep=3pt
      \begin{array}{ccccccccccc}
        0 & 0 & 0 & 0 & 0 & 48510 & 95223 & 0 & 0 & 0 & 16170 \\
        0 & 0 & 0 & 425810 & 0 & 2817045 & 9741270 & 1348462 & 0 & 0 & 1252020 \\
        0 & 0 & 3476550 & 6110115 & 48434732 & 0 & 336258384 & 218861747 & 0 & 0 & 38848656 \\
        0 & 6055665 & 95465370 & 190273710 & 1221447150 & 1171139970 & 2726237052 & 7298343195 & 0 & 0 & 1158647028 \\
        3128160 & 204553195 & 1480047030 & 5336244870 & 15771654360 & 32618391960 & 21628131756 & 73442566505 & 29020437724 & 0 & 0 \\
        116263280 & 3318028175 & 18873326010 & 77596724865 & 174458095350 & 356484794820 & 298526146072 & 392659895320 & 419332019003 & 0 & 0 \\
        1988081620 & 35746404925 & 197029544250 & 726747795015 & 1551387171180 & 2562636437130 & 2754345379016 & 1907270574440 & 2403530867430 & 384166681454 & 0 \\
        21102099620 & 290306961329 & 1558206940290 & 4960832042100 & 10157416978170 & 14397653320512 & 15830016304564 & 10335346465675 & 8454107203080 & 3721378398370 & 0 \\
        157521933940 & 1842856573327 & 9151953918030 & 25828662738780 & 49078270584420 & 64463871381756 & 64790433176084 & 46215397662665 & 25413887653770 & 13866272170542 & 3624993260826 \\
        879576036500 & 9097329993521 & 40326294432270 & 103410904320900 & 179681190528840 & 223720508502183 & 206003553429058 & 148200002432020 & 74515561079190 & 38329760467746 & 21352330210512 \\
        3768921407020 & 34425402760287 & 134937782918400 & 317406347163180 & 506548245985320 & 592781349468231 & 516092578758680 & 349390783269990 & 182484813042840 & 86130722275092 & 36302262824520 \\
        12408373123020 & 98585568531450 & 344523626901300 & 742720332283380 & 1098864688687920 & 1195101212250330 & 988299659161584 & 622626376181040 & 315159824447160 & 127817022168000 & 27594339093216 \\
        30888195414300 & 211221314490186 & 667161153364860 & 1314049020225480 & 1805851509808500 & 1815258345062778 & 1389893720693940 & 808485195464100 & 350770360746960 & 103510976206176 & 7959911420160 \\
        56418683248620 & 333956069661060 & 961433219937960 & 1730658737031840 & 2184547015159740 & 2011470759046980 & 1371614133582000 & 691021013177880 & 227021138467200 & 33927611477760 & 0 \\
        72303376012560 & 380597496158880 & 996948009953280 & 1640782277386560 & 1861980233153040 & 1520944502505120 & 876545356468320 & 330199345808640 & 64746646272000 & 0 & 0 \\
        60744201708960 & 297669661581600 & 704885552078400 & 1045548136987200 & 1036553459911200 & 695810390758560 & 300533746867200 & 63688771238400 & 0 & 0 & 0 \\
        29689237670400 & 143165195712000 & 307454361984000 & 391950244224000 & 317437952832000 & 151152068390400 & 31685955840000 & 0 & 0 & 0 & 0 \\
        6337191168000 & 31685955840000 & 63371911680000 & 63371911680000 & 31685955840000 & 6337191168000 & 0 & 0 & 0 & 0 & 0 \\
      \end{array}
      \right) $
  \end{Exx} }

\medskip\noindent
\resizebox{\linewidth}{!}{
  \begin{Exx}
    {\large $ \text{M}^{\rm T} = $}
    $ \left( \arraycolsep=3pt
      \begin{array}{ccccccccccc}
        0 & 0 & 0 & 0 & 0 & 48510 & 97020 & 0 & -64680 & 0 & 16170 \\
        0 & 0 & 0 & 544005 & -709170 & 2817045 & 9741270 & 2279970 & -5453910 & -810810 & 1252020 \\
        0 & 0 & 3476550 & 6110115 & 51415980 & -17887485 & 336258384 & 259149660 & -233284590 & -50657310 & 38848656 \\
        0 & 6055665 & 95465370 & 190273710 & 1221447150 & 1171139970 & 2726237052 & 7522825695 & -627553080 & -4316051520 & 1158647028 \\
        3128160 & 204553195 & 1480047030 & 5336244870 & 15771654360 & 32618391960 & 21628131756 & 73442566505 & 36182631870 & -42053593230 & -5517429876 \\
        116263280 & 3318028175 & 18873326010 & 77596724865 & 174458095350 & 356484794820 & 298526146072 & 392659895320 & 439184120760 & -87653879280 & -188752387572 \\
        1988081620 & 35746404925 & 197029544250 & 726747795015 & 1551387171180 & 2562636437130 & 2754345379016 & 1907270574440 & 2403530867430 & 568697131530 & -1107182700456 \\
        21102099620 & 290306961329 & 1558206940290 & 4960832042100 & 10157416978170 & 14397653320512 & 15830016304564 & 10335346465675 & 8454107203080 & 4052140160904 & -1984570575204 \\
        157521933940 & 1842856573327 & 9151953918030 & 25828662738780 & 49078270584420 & 64463871381756 & 64790433176084 & 46215397662665 & 25413887653770 & 13866272170542 & 3624993260826 \\
        879576036500 & 9097329993521 & 40326294432270 & 103410904320900 & 179681190528840 & 223720508502183 & 206003553429058 & 148200002432020 & 74515561079190 & 38329760467746 & 21352330210512 \\
        3768921407020 & 34425402760287 & 134937782918400 & 317406347163180 & 506548245985320 & 592781349468231 & 516092578758680 & 349390783269990 & 182484813042840 & 86130722275092 & 36302262824520 \\
        12408373123020 & 98585568531450 & 344523626901300 & 742720332283380 & 1098864688687920 & 1195101212250330 & 988299659161584 & 622626376181040 & 315159824447160 & 127817022168000 & 27594339093216 \\
        30888195414300 & 211221314490186 & 667161153364860 & 1314049020225480 & 1805851509808500 & 1815258345062778 & 1389893720693940 & 808485195464100 & 350770360746960 & 103510976206176 & 7959911420160 \\
        56418683248620 & 333956069661060 & 961433219937960 & 1730658737031840 & 2184547015159740 & 2011470759046980 & 1371614133582000 & 691021013177880 & 227021138467200 & 33927611477760 & 0 \\
        72303376012560 & 380597496158880 & 996948009953280 & 1640782277386560 & 1861980233153040 & 1520944502505120 & 876545356468320 & 330199345808640 & 64746646272000 & 0 & 0 \\
        60744201708960 & 297669661581600 & 704885552078400 & 1045548136987200 & 1036553459911200 & 695810390758560 & 300533746867200 & 63688771238400 & 0 & 0 & 0 \\
        29689237670400 & 143165195712000 & 307454361984000 & 391950244224000 & 317437952832000 & 151152068390400 & 31685955840000 & 0 & 0 & 0 & 0 \\
        6337191168000 & 31685955840000 & 63371911680000 & 63371911680000 & 31685955840000 & 6337191168000 & 0 & 0 & 0 & 0 & 0 \\
      \end{array}
      \right) $
 \end{Exx} }


\subsection{Upper bounds for \boldmath{$c_n(\alpha)$}}

We apply Proposition~\ref{p2.3} to estimate the largest zero
$\,x_n=c_n^2(\al)\,$ of the polynomial $\,R_n(x)\,$ from above,
$$
  x_n\leq u_k(R_n)=p_k(R_n)^{1/k}\,,\qquad k=3,\,4,\,5,\,6\,.
$$
Then with the assistance of computer algebra we obtain a further
estimation of the form
$$
  u_k(R_n)\leq c^{1/k}\,(n+1)(n+\sigma(\alpha+1)),
$$
with the optimal (i.e., the smallest possible) constants
$\,c=c(k)\,$ and $\,\sigma=\sigma(k)$.

The algorithm is analogous to Algorithm~\ref{alg1}, and the code has
only a few differences which are specified later.

\begin{algorithm}
\caption{\quad Estimating $\,c_n(\alpha)\,$ from above} \label{alg2}
\begin{tabular}{ll}
  {\sl Input:}  & $\,k\in\{3,4,5,6\}\,$ -- the number of the highest degree coefficients of $\,R_n(x)\,$ \\
  {\sl Step 1.} & Express the power sum $\,p_k(R_n)\,$ in terms of $\,\{b_i\}_{i=1}^k\,$ \\
  {\sl Step 2.} & Find $\,\{b_i\}_{i=1}^k\,$ in terms of $\,n\,$ and $\,\alpha\,$ using Proposition~\ref{p2.2} \\
  {\sl Step 3.} & Find a proper value $\sigma$ for parameter $\,s\,$ in the expression \\
                & $\,c\,(n+1)^k(n+s(\alpha+1))^k-p_k\,$, \ where $\,c\,$ is the coefficient of $\,n^{2k}\,$ in $\,p_k\,$ \\
  {\sl Step 4.} & Represent the numerator of $\,f=c\,(n+1)^k(n+\sigma(\alpha+1))^k-p_k\,$ \\
                & in powers of $\,n\,$ and $\,(\alpha+1)\,$ \\
  {\sl Step 5.} & Estimate from below the expression $\,f\,$ to prove that $\,f\ge 0\,$
\end{tabular}
\end{algorithm}

\smallskip
\underline{\sl Step 1:} \ The same as in Algorithm~\ref{alg1}.

\smallskip
\underline{\sl Step 2:} \ Identical to that in Algorithm~\ref{alg1}.

\smallskip
\underline{\sl Step 3:} \ The only differences with Algorithm~\ref{alg1}
are that we set $\,c\,$ to be the coefficient of $\,n^{2k}\,$ in $\,p_k\,$ and
$$
  f_s=c\,(n+1)^k(n+s(\alpha+1))^k-p_k\,.
$$

The highest order coefficients in $\,\sum_j
\mu_{i,j}(s)(\alpha+1)^{d-j}\,$ are functions in $\,s\,$ of the form
$\,A_i s^\nu-B_i\,$, with $\,A_i>0\,$ and $\,B_i\ge 0\,$. We denote
their non-negative zeros by $\,s_i\,$ for each $\,i\,$ and choose
$\,\sigma=\max_i s_i\,$.

The results for $\,k=3,4,5,6\,$ obtained by symbolic computations
are given in Table~\ref{tab2}. \vspace{-2ex}
\begin{table}[H]
  \hspace*{5em}\caption{The optimal values of $\,c\,$ and $\,\sigma\,$ in the upper bounds for $c_n^2(\al)$.} \label{tab2}
  \begin{center}
  \begin{tabular}{ccc} \hline
      $k$ \qquad & $c$ & \qquad $\sigma$ \\[3pt] \hline
       3  \qquad & $\dfrac{1}{(\alpha+1)^3(\alpha+3)(\alpha+5)}$ & \qquad \mbox{\rule[-4mm]{0mm}{10mm}{$\dfrac25$}} \\[8pt]
       4  \qquad & $\dfrac{5\alpha+17}{2(\alpha+1)^4(\alpha+3)^2(\alpha+5)(\alpha+7)}$ & \qquad $\dfrac37$ \\[8pt]
       5  \qquad & $\dfrac{(7\alpha+31)}{(\alpha+1)^5(\alpha+3)^2(\alpha+5)(\alpha+7)(\alpha+9)}$ & \qquad $\dfrac49$ \\[8pt]
       6  \qquad & \quad $\dfrac{21\alpha^3+299\alpha^2+1391\alpha+2073}{(\alpha+1)^6(\alpha+3)^3(\alpha+5)^2(\alpha+7)(\alpha+9)(\alpha+11)}$
                     & \qquad $\dfrac{5}{11}$ \\[8pt] \hline
  \end{tabular}
    \end{center}
\end{table}

\underline{\sl Step 4:} \
With $\,c\,$ and $\,\sigma\,$ determined in the previous Step~3 we set
$$
  f=c\,(n+1)^k(n+\sigma(\alpha+1))^k-p_k=:\frac{\varphi(n,\alpha)}{\psi(\alpha)}\,.
$$
The rest of the source has no difference with Step~4 of Algorithm~\ref{alg1}.

\smallskip
\underline{\sl Step 5:} \
The same as in Algorithm~\ref{alg1}. Using the same recursive procedure
we find a matrix  $\,\mathbf{\Lambda}\,$ satisfying
$\,\mathbf{0}\le\mathbf{\Lambda}\le\mathbf{M}\,$. Then
$\,\varphi(n,\alpha)\ge 0\,$, $\,f\ge 0\,$ and therefore
$$
  c_n^{2k}(\alpha) \le p_k \le c\,(n+1)^k(n+\sigma(\alpha+1))^k\,
$$
for the corresponding $\,c\,$ and $\,\sigma\,$ evaluated in Step~3.
For $\,k=3,4,5,6\,$ we obtain estimations \eqref{e1.10}, \eqref{e1.12},
\eqref{e1.14}, and \eqref{e1.16}, respectively.

The matrices $\,\mathbf{M}\,$ from Step~4 and $\,\mathbf{\Lambda}\,$
from Step~5 obtained with  {\sl Mathematica} are given below.

\bigskip
Case $\,k=3\,$:

\begin{center}
  \begin{Exx}
    {\scriptsize
    $ \mathbf{\Lambda} =
      \left( \arraycolsep=3pt
      \begin{array}{ccccc}
        0 & 0 & 0 & 1500 & 3300 \\
        0 & 115 & 1885 & 4170 & 4233 \\
        32 & 598 & 3026 & 6360 & 0 \\
        96 & 979 & 2143 & 850 & 0 \\
        96 & 624 & 1098 & 0 & 0 \\
      \end{array}
      \right)$
    \qquad\qquad
    $ \text{M} =
      \left( \arraycolsep=3pt
      \begin{array}{ccccc}
        0 & 0 & 0 & 1500 & 3300 \\
        0 & 115 & 1885 & 4170 & 4650 \\
        32 & 598 & 3026 & 6360 & -600 \\
        96 & 979 & 2143 & 1560 & -1950 \\
        96 & 624 & 1098 & -2130 & 0 \\
      \end{array}
    \right) $}
  \end{Exx}
\end{center}

\bigskip
Case $\,k=4\,$:

\medskip
\begin{Exx}
  {\scriptsize
  $ \mathbf{\Lambda} =
    \left( \arraycolsep=3pt
    \begin{array}{ccccccccc}
      0 & 0 & 0 & 0 & 905520 & 8808240 & 29717520 & 41571600 & 19756800 \\
      0 & 0 & 54390 & 2038890 & 16676660 & 60285680 & 115770830 & 117031110 & 48774600 \\
      0 & 42294 & 1237572 & 10966494 & 52723608 & 141477042 & 198565500 & 127823850 & 24194362 \\
      6075 & 266115 & 3694950 & 25364010 & 85166735 & 157047575 & 154257320 & 46893642 & 0 \\
      24300 & 617510 & 5700800 & 26734470 & 72437020 & 97039330 & 34815501 & 0 & 0 \\
      36450 & 678780 & 4979940 & 16392810 & 28823750 & 17907835 & 0 & 0 & 0 \\
      24300 & 360421 & 2131108 & 6792156 & 5246162 & 0 & 0 & 0 & 0 \\
    \end{array}
    \right) $

  \medskip
  $ \text{M} =
    \left( \arraycolsep=3pt
    \begin{array}{ccccccccc}
      0 & 0 & 0 & 0 & 905520 & 8808240 & 29717520 & 41571600 & 19756800 \\
      0 & 0 & 54390 & 2038890 & 16676660 & 60285680 & 115770830 & 117031110 & 48774600 \\
      0 & 42294 & 1237572 & 10966494 & 52723608 & 141477042 & 198565500 & 127823850 & 27783000 \\
      6075 & 266115 & 3694950 & 25364010 & 85166735 & 157047575 & 154257320 & 52558380 & -11730600 \\
      24300 & 617510 & 5700800 & 26734470 & 72437020 & 97039330 & 38636640 & -18088350 & -10495800 \\
      36450 & 678780 & 4979940 & 16392810 & 28823750 & 20280800 & -12849340 & -18282390 & 0 \\
      24300 & 360421 & 2131108 & 6792156 & 5246162 & -9491857 & -9740850 & 0 & 0 \\
    \end{array}
    \right) $ }
\end{Exx}

\bigskip
Case $\,k=5\,$:

\medskip\noindent
\scalebox{0.6}{ \begin{Exx}
                  $ \mathbf{\Lambda} =
                    \left( \arraycolsep=3pt
                    \begin{array}{ccccccccccc}
                      0 & 0 & 0 & 0 & 0 & 85424220 & 1436596560 & 8988832440 & 26097558480 & 34662943980 & 16203045600 \\
                      0 & 0 & 0 & 4261005 & 260814330 & 3617057430 & 22250151630 & 73071107235 & 134891273160 & 134642808090 & 56710659600 \\
                      0 & 0 & 5436720 & 241567920 & 3235204800 & 22246774740 & 91003127400 & 223063050420 & 312360753600 & 222393230640 & 64812182400 \\
                      0 & 1982358 & 88937982 & 1392482448 & 12340605438 & 63755213760 & 194677526736 & 357163148790 & 375802372260 & 186521488020 & 12638375568 \\
                      200704 & 14563010 & 340432890 & 4020858058 & 25446365294 & 99455228208 & 241336266948 & 338611016520 & 235926284580 & 44541786567 & 0 \\
                      1003520 & 42390775 & 693405300 & 6004806185 & 31876009900 & 96870254355 & 175080003840 & 176585507595 & 54286938720 & 0 & 0 \\
                      2007040 & 63580160 & 829630410 & 5638883530 & 22495811450 & 57112266330 & 77686343280 & 30853075478 & 0 & 0 & 0 \\
                      2007040 & 52428341 & 568553244 & 3375204826 & 9950248616 & 17535199185 & 13032227178 & 0 & 0 & 0 & 0 \\
                      1003520 & 22758400 & 207566490 & 998218460 & 3486984100 & 3092469120 & 0 & 0 & 0 & 0 & 0 \\
                    \end{array}
                    \right) $
                \end{Exx} }

\medskip\noindent
\scalebox{0.6}{ \begin{Exx}
                  $ \text{M} =
                    \left( \arraycolsep=3pt
                    \begin{array}{ccccccccccc}
                      0 & 0 & 0 & 0 & 0 & 85424220 & 1436596560 & 8988832440 & 26097558480 & 34662943980 & 16203045600 \\
                      0 & 0 & 0 & 4261005 & 260814330 & 3617057430 & 22250151630 & 73071107235 & 134891273160 & 134642808090 & 56710659600 \\
                      0 & 0 & 5436720 & 241567920 & 3235204800 & 22246774740 & 91003127400 & 223063050420 & 312360753600 & 222393230640 & 64812182400 \\
                      0 & 1982358 & 88937982 & 1392482448 & 12340605438 & 63755213760 & 194677526736 & 357163148790 & 375802372260 & 186521488020 & 16203045600 \\
                      200704 & 14563010 & 340432890 & 4020858058 & 25446365294 & 99455228208 & 241336266948 & 338611016520 & 235926284580 & 51689001420 & -16203045600 \\
                      1003520 & 42390775 & 693405300 & 6004806185 & 31876009900 & 96870254355 & 175080003840 & 176585507595 & 59214803760 & -31849915230 & -8101522800 \\
                      2007040 & 63580160 & 829630410 & 5638883530 & 22495811450 & 57112266330 & 77686343280 & 32878980540 & -21278795580 & -19430795160 & 0 \\
                      2007040 & 52428341 & 568553244 & 3375204826 & 9950248616 & 17535199185 & 14090589072 & -8987585040 & -16802648100 & 0 & 0 \\
                      1003520 & 22758400 & 207566490 & 998218460 & 3486984100 & 3092469120 & -5291809470 & -5709701340 & 0 & 0 & 0 \\
                    \end{array}
                    \right) $
                \end{Exx} }

\bigskip
Case $\,k=6\,$:

\medskip\noindent
\resizebox{\linewidth}{!}{
  \begin{Exx}
    {\large $ \mathbf{\Lambda^T} = $}
    $ \left( \arraycolsep=3pt
      \begin{array}{ccccccccccc}
        0 & 0 & 0 & 0 & 0 & 137812500 & 826875000 & 2067187500 & 2756250000 & 2067187500 & 826875000 \\
        0 & 0 & 0 & 0 & 1712831340 & 16225847190 & 61013597190 & 121629375000 & 141073384270 & 96161625000 & 35900159390 \\
        0 & 0 & 0 & 7214978925 & 129414811020 & 682431963570 & 1909641192060 & 3176038395495 & 3233829060980 & 1990746645930 & 685515044860 \\
        0 & 0 & 11987252200 & 496609292025 & 3975634217280 & 15545520413640 & 34393961720958 & 48625624140345 & 44026671755710 & 24295459922370 & 7604368197232 \\
        0 & 6757992780 & 889711262440 & 13291579355670 & 71880656791680 & 218617776740580 & 407202449242428 & 488239798537050 & 391807489052700 & 198688598846220 & 53934381373052 \\
        0 & 678222070680 & 23583750724380 & 204179315141340 & 848078318563740 & 2083353346583670 & 3322865459747106 & 3448886831378940 & 2400975627767010 & 1108865878357650 & 280425349855044 \\
        156657528720 & 19589397975840 & 339201489248160 & 2064129233162535 & 6887056367079900 & 14142830057025810 & 19301411805553332 & 17606432763060855 & 10537484120843400 & 4243406630751540 & 1099649486098068 \\
        5737886375760 & 285634042298220 & 3098176511479280 & 14556440238205635 & 39911599747186200 & 69897865311080520 & 81762211401327198 & 65213552886377595 & 33832857309375390 & 11331273200216130 & 2850997375559272 \\
        90962828787600 & 2522414547688728 & 19367097472935520 & 73609280525046420 & 168732748048142520 & 253532902195144824 & 255297208204811052 & 174944618786085540 & 78706656458079300 & 21128915605090608 & 3776309107237628 \\
        821967442647120 & 14663875740671388 & 85966291726204120 & 270402692712289830 & 525000679384868280 & 675267348901207524 & 584508472246139364 & 337073370501014070 & 126378098307624240 & 25176310282435158 & 100310771384846 \\
        4683403445822640 & 58501440227708016 & 275002627721514000 & 724008553578199890 & 1198318687198670400 & 1310711831601226488 & 960008914211603088 & 452338940602722600 & 125459625652811660 & 10624247696162964 & 0 \\
        17541804701337840 & 163077338127380292 & 633513827175290600 & 1404602858405626110 & 1976258187012861840 & 1815905828323340796 & 1083067213944785244 & 386368902093765990 & 50117292233704841 & 0 & 0 \\
        43597315400007600 & 317712831386568156 & 1035759237211329400 & 1939380788446101000 & 2286621788179676640 & 1716811248944619828 & 768783125513892720 & 147265658296879789 & 0 & 0 & 0 \\
        70706025735594480 & 424476237961133820 & 1166286653916203820 & 1837004374966081860 & 1763385304898752920 & 1006399215292377060 & 256729802815313776 & 0 & 0 & 0 & 0 \\
        71282695085965440 & 370662526010533680 & 857715825988763280 & 1112754260403129960 & 820925242403999040 & 277604251541734810 & 0 & 0 & 0 & 0 & 0 \\
        40151827863100800 & 190176873072832800 & 371893024262944800 & 374764348674777600 & 176396296906922400 & 0 & 0 & 0 & 0 & 0 & 0 \\
        9572830210944000 & 43261828837920000 & 72716691025440000 & 49617294906564000 & 0 & 0 & 0 & 0 & 0 & 0 & 0 \\
      \end{array}
      \right) $
  \end{Exx} }

\medskip\noindent
\resizebox{\linewidth}{!}{
  \begin{Exx}
    {\large $ \text{M}^{\rm T} = $}
    $ \left( \arraycolsep=3pt
      \begin{array}{ccccccccccc}
        0 & 0 & 0 & 0 & 0 & 137812500 & 826875000 & 2067187500 & 2756250000 & 2067187500 & 826875000 \\
        0 & 0 & 0 & 0 & 1712831340 & 16225847190 & 61013597190 & 121629375000 & 141073384270 & 96161625000 & 35900159390 \\
        0 & 0 & 0 & 7214978925 & 129414811020 & 682431963570 & 1909641192060 & 3176038395495 & 3233829060980 & 1990746645930 & 685515044860 \\
        0 & 0 & 11987252200 & 496609292025 & 3975634217280 & 15545520413640 & 34393961720958 & 48625624140345 & 44026671755710 & 24295459922370 & 7604368197232 \\
        0 & 6757992780 & 889711262440 & 13291579355670 & 71880656791680 & 218617776740580 & 407202449242428 & 488239798537050 & 391807489052700 & 198688598846220 &  53934381373052 \\
        0 & 678222070680 & 23583750724380 & 204179315141340 & 848078318563740 & 2083353346583670 & 3322865459747106 & 3448886831378940 & 2400975627767010 & 1108865878357650 & 280425349855044 \\
        156657528720 & 19589397975840 & 339201489248160 & 2064129233162535 & 6887056367079900 & 14142830057025810 & 19301411805553332 & 17606432763060855 & 10537484120843400 & 4243406630751540 & 1099649486098068 \\
        5737886375760 & 285634042298220 & 3098176511479280 & 14556440238205635 & 39911599747186200 & 69897865311080520 & 81762211401327198 & 65213552886377595 & 33832857309375390 & 11331273200216130 & 2850997375559272 \\
        90962828787600 & 2522414547688728 & 19367097472935520 & 73609280525046420 & 168732748048142520 & 253532902195144824 & 255297208204811052 & 174944618786085540 & 78706656458079300 & 21128915605090608 & 3776309107237628 \\
        821967442647120 & 14663875740671388 & 85966291726204120 & 270402692712289830 & 525000679384868280 & 675267348901207524 & 584508472246139364 & 337073370501014070 & 126378098307624240 & 25176310282435158 & 100310771384846 \\
        4683403445822640 & 58501440227708016 & 275002627721514000 & 724008553578199890 & 1198318687198670400 & 1310711831601226488 & 960008914211603088 & 452338940602722600 & 125459625652811660 & 11652636597863226 & -6170333410201568 \\
        17541804701337840 & 163077338127380292 & 633513827175290600 & 1404602858405626110 & 1976258187012861840 & 1815905828323340796 & 1083067213944785244 & 386368902093765990 & 52601940894289780 & -13776737928943668 & -6786924207395784 \\
        43597315400007600 & 317712831386568156 & 1035759237211329400 & 1939380788446101000 & 2286621788179676640 & 1716811248944619828 & 768783125513892720 & 152704407947442120 & -28750678275282040 & -22906285043773824 & -2307796348667040 \\
        70706025735594480 & 424476237961133820 & 1166286653916203820 & 1837004374966081860 & 1763385304898752920 & 1006399215292377060 & 264555646816493340 & -39077291676501540 & -45675842674954800 & -9544747851001440 & 0 \\
        71282695085965440 & 370662526010533680 & 857715825988763280 & 1112754260403129960 & 820925242403999040 & 283521195216151200 & -25185070761515280 & -59085225024570360 & -16885936111968000 & 0 & 0 \\
        40151827863100800 & 190176873072832800 & 371893024262944800 & 374764348674777600 & 178350638317790400 & -4068645433221600 & -43598489040136800 & -14075574910689600 & 0 & 0 & 0 \\
        9572830210944000 & 43261828837920000 & 72716691025440000 & 49705079941440000 & 1840928886720000 & -13438780873056000 & -4602322216800000 & 0 & 0 & 0 & 0 \\
      \end{array}
      \right)$
  \end{Exx} }


\bigskip
\setcounter{equation}{0}
\section{Concluding remarks}

\textbf{1.} In our computer algebra approach for derivation of
bounds for the best Markov constant $\,c_n(\al)\,$ we perform some
optimization with respect to  parameter~$\,s\,$. Our motivation for
searching lower bounds for $\,c_n^2(\al)\,$ with a factor depending
on $\,n\,$ of the special form $\,n\big(n+\sigma(\al+1)\big)\,$ is
Corollary~D(ii).

An interesting observation about the lower bounds
$\,\underline{c}_{n,k}(\al)\,$ in Theorem~\ref{t1.1} is that they imply
$$
  \frac{k\,n}{k+1} =
  \lim_{\al\to \infty}\al\,\underline{c}_{n,k}^{2}(\al) \le
  \lim_{\al\to \infty}\al\,c_n^{2}(\al)\,,
  \qquad 3\leq k\leq 6
$$
(the lower bound in Corollary~\ref{c1.4} follows from the case
$\,k=6\,$). This observation and Proposition~\ref{p2.3} give rise
for the following

\begin{conjecture}\label{con4.1}
The best Markov constant $\,c_n(\al)\,$ satisfies the asymptotic relation:
$$
  \lim_{\al\to \infty}\al\,c_n^{2}(\al) = n\,.
$$
\end{conjecture}

We also performed a search for lower bounds for $\,c_n^2(\al)\,$
with a factor depending on $\,n\,$ of the form
$\,(n+1)\big(n+\sigma(\al+1)\big)\,$. Such a choice is reasonable,
as the resulting lower bounds preserve the limit relation in
Corollary~D\,(i). The optimal value then is $\,\sigma=-1/3\,$ (the
same for all $\,k\,$, $\,3\leq k\leq 6$), and we obtain lower bounds
as in Theorem~\ref{t1.1} with $\,n\big(n+\sigma(\al+1)\big)\,$
replaced by $\,(n+1)\big(n-(\al+1)/3\big)\,$\,. These lower bounds
make sense only for $n>(\al+1)/3$, and are better than those in
Theorem~\ref{t1.1} only for $\,\al\,$ close to $\,-1\,$.

\medskip
\textbf{2.} The bounds $\,\big(\underline{c}_{\,n,k}(\alpha),
\overline{c}_{\,n,k}(\alpha)\big)\,$ ($\,3\leq k \leq 6$) in
Theorem~\ref{t1.1} imply bounds $\,(\ell_k(\al),u_k(\al))\,$
(occurring in the middle columns of Tables~1 and 2) for the
asymptotic Markov constant $\,c(\al)\,$, and the bounds deduced with
a larger $\,k\,$ are superior. While the lower bounds
$\,\ell_k(\al)\,$ are of the correct order $\,\OO(\al^{-1})\,$ as
$\,\al\to\infty$, for the upper bound $\,u_k(\al)\,$ we have
$\,u_k(\al)=\OO(\al^{-1+\frac{1}{2k}})\,$ as $\,\al\to\infty$,
$\,(3\leq k\leq 6)\,$. The ratio
$$
  \rho_k(\al) := \frac{u_k(\al)}{\ell_k(\al)}\,,
    \qquad 3\leq k\leq 6,
$$
tends to $\,1\,$ as $\,\al\to -1$, which indicates that for moderate
$\,\al\,$ the bounds $\,\ell_k(\al)\,$ and $\,u_k(\al)\,$ are rather
tight. This observation is clearly seen in the particular case
$\,\al=0\,$, where, according to Tur\'{a}n's result, we have
$\,c(0)=\frac{2}{\pi}\,$. We give the lower and the upper bounds for
$\,c(0)\,$ and the overestimation factors in Table~\ref{tab3}.
\vspace{-3ex}
\begin{table}[H]
  \hspace*{5em}\caption{The lower and the upper bounds for the asymptotic
  Markov constant $\,c(0)\,$ and the overestimation factors.} \label{tab3}
  \begin{center}
  \begin{tabular}{ccccc} \hline
    $k$ & \!\! $\ell_k(0)$ \!\! & \!\! $u_k(0)$ \!\! & \!\! $\dfrac{c(0)}{\ell_k(0)}$ \!\!
                   & \!\! $\dfrac{u_k(0)}{c(0)}$ \!\! \\[8pt] \hline
     3  & \!\! $\sqrt{\frac{2}{5}}\approx 0.63245553$ \!\! & \!\! $\sqrt[6]{\frac{1}{15}}\approx 0.63677321$ \!\!
                   & \!\! $1.006584242$ \!\! & $1.00024103$ \!\! \\[8pt]
     4  & \!\! $\sqrt{\frac{17}{42}}\approx 0.63620901$ \!\! & $\sqrt[8]{\frac{17}{630}}\approx 0.63663212$ \!\!
                   & \!\! $1.00064564$ \!\! & $1.00001939$ \!\! \\[8pt]
     5  & \!\! $\sqrt{\frac{62}{153}}\approx 0.63657580$ \!\! & \!\! $\sqrt[10]{\frac{31}{2835}}\approx 0.63662085$ \!\!
                   & \!\! $1.00006906$\!\! & $1.00000170$ \!\! \\[8pt]
     6  & \!\! $\sqrt{\frac{2073}{5115}}\approx 0.63661494$ \!\! & \!\! $\sqrt[12]{\frac{2073}{467775}}\approx 0.63661987$ \!\!
                   & \!\! $1.00000757$ \!\! & $1.00000015$ \!\! \\[8pt]
    \hline
  \end{tabular}
    \end{center}
\end{table}

Although the ratios $\,\rho_k\,$, $\,3\leq k\leq 6\,$, satisfy
$\,\rho_k(\al)\to\infty\,$ as $\,\al\to\infty\,$, they grow rather
slowly. For instance, $\,\rho_6(\al)<2\,$ for $\,\al<140000\,$, see
Figure 1.

\begin{figure}[htp]
  \centering
  \includegraphics[scale=0.85,clip]{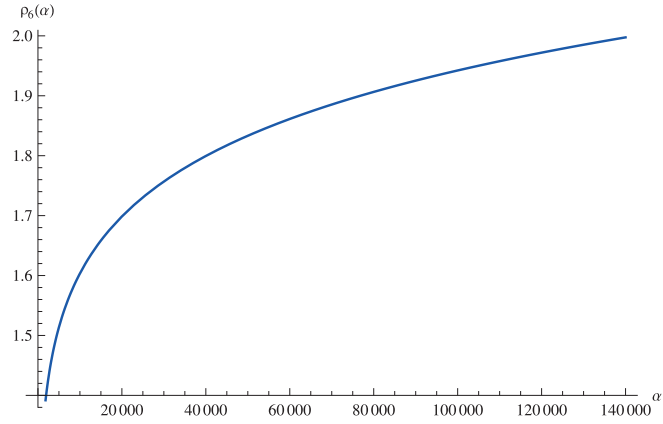}\hfill
  \caption{The graph of $\,\rho_6(\al)<2\,$.} \label{fig1}
\end{figure}

\pagebreak
\textbf{3.} Another interesting observation, concerning the coefficients
of $\,R_n\,$ inspires the following

\begin{conjecture} \label{con4.2}
For every fixed $\,k\in \mathbb{N}\,$, the coefficient $\,b_{k,n}\,$,
$\,n>k\,$, of the polynomial
$$
  \,R_n(x)=x^n-b_{1,n}\,x^{n-1}+b_{2,n}\,x^{n-2}-\cdots+(-1)^n\,b_{n,n}\,,
$$
satisfies
\begin{equation} \label{e4.1}
  b_{k,n} = \frac{n^{2k}}{2^k\,k!(\al+1)\cdots(\al+2k-1)}+\OO(n^{2k-1})\,.
\end{equation}
\end{conjecture}

Conjecture~\ref{con4.2} is verified with our computer algebra
approach for $\,1\leq k\leq 6\,$, but so far we do not have a proof
for the general case. Having \eqref{e4.1} proved, we could try to
find the explicit form of $\,d_k\,$, the coefficient of $\,n^{2k}\,$
in Newton's function $\,p_k(R_n)\,$, and consequently to obtain two
sequences $\{\ell_k\}\,$ and $\,\{u_k\}\,$ defined by
$\,\ell_k=\sqrt{d_k/d_{k-1}}\,$ and $\,u_k=\sqrt[2k]{d_k}\,$ which tend
monotonically from below and from above, respectively, to
$\,c(\al)\,$, the sharp asymptotic Markov constant.


\bigskip
{\bf Acknowledgement.} The authors are supported by the
Bulgarian National Research Fund under Contract DN 02/14 and by the
Sofia University Research Fund under Contract 80.10-11/2017.



\bigskip\bigskip\bigskip
\noindent
{\sc Geno Nikolov, \ \ Rumen Uluchev} \smallskip\\
Department of Mathematics and Informatics\\
University of Sofia \\
5 James Bourchier Blvd. \\
1164 Sofia \\
BULGARIA \\
{\it E-mails:} \ {\tt geno@fmi.uni-sofia.bg}, \ \ {\tt rumenu@fmi.uni-sofia.bg}


\end{document}